\documentclass{TEMA}

\usepackage[brazil]{babel}      
\usepackage[utf8]{inputenc}    

\usepackage[dvips]{graphics}
\usepackage{subfigure}
\usepackage[pdftex]{graphicx}
\usepackage{epsfig}
\usepackage{hyperref}
\usepackage{framed}
\usepackage{psfrag}
\usepackage{tikz} 

\usepackage{mathrsfs}
\usepackage{amssymb}

\begin{document}
\title{Método de Monte Carlo aplicado ao Cálculo Fracionário}

\author {
 L.N. FERREIRA%
 \thanks{luverci@gmail.com}, 
 Instituto de Matemática, Estatística e Física,
 IMEF, Universidade Federal do Rio Grande,
 Av. Itália, km 8, bairro Carreiros, 
96203-900, Rio Grande, RS, Brasil.
\\ \\

M.J. LAZO%
\thanks{matheusjlazo@gmail.com}, 
 Instituto de Matemática, Estatística e Física,
 IMEF, Universidade Federal do Rio Grande,
 Av. Itália, km 8, bairro Carreiros, 
96203-900, Rio Grande, RS, Brasil.
}
\criartitulo
\runningheads{L.N. Ferreira e M.J. Lazo}
{Método de Monte Carlo aplicado ao Cálculo Fracionário}

\begin{abstract}
{\bf Resumo} O presente trabalho analisa e desenvolve um método para resolver equações diferenciais fracionárias utilizando o Método de Monte Carlo. Uma simulação numérica é realizada para algumas equações diferenciais, comparando os resultados com o que existe na literatura matemática. A linguagem Python é usada para criar modelos computacionais.

{\bf Palavras-chave:} Cálculo Fracionário, Monte Carlo, Python.
\end{abstract}

  \section{Introdução}

  O Cálculo Fracionário, como são chamadas as generalizações do cálculo para derivadas e integrais de ordem não inteira, teve seu início em 1695 quando Leibniz escreveu uma carta à l'H\^opital propondo uma derivada de ordem $\frac{1}{2}$ \cite{OldhamSpanier,Podlubny}. Outros grandes matemáticos da história também contribuíram para o desenvolvimento do Cálculo Fracionário, entre eles podemos citar nomes como Euler, Laplace, Liouville, Gr$\ddot{\mbox{u}}$nwald, Letnikov, Riemann e muitos outros. Apesar de ser quase tão antigo quanto o Cálculo usual de derivadas e integrais de ordem inteira, só nas últimas quatro décadas o Cálculo Fracionário despertou maior atenção devido a suas aplicações em vários campos da ciência e da engenharia. Esta demora no surgimento de aplicações do Cálculo Fracionário se explica devido as dificuldades de se interpretar física e geometricamente as derivadas e integrais não inteiras. Até hoje, mais de trezentos anos após sua criação, ainda não temos uma interpretação física e uma interpretação geométrica consistentes para essas derivadas e integrais. Recentemente foi proposta uma interpretação geométrica como ``sombras de uma área'' e uma interpretação física como ``sombras do passado'' \cite{Podlubny}. Embora essas interpretações ainda pareçam não totalmente satisfatórias, ficou evidente nestas últimas quatro décadas que o Cálculo Fracionário é uma ferramenta matemática extremamente importante para aplicações que envolvam o estudo de sistemas complexos \cite{mainardi, Hilfer, SATM}. 
  
  As derivadas fracionárias são, em geral, operadores não-locais. Consequentemente, os sistemas complexos em que utilizamos o Cálculo Fracionário são sistemas com comportamento não-local ou com dinâmica dependente de memória. Entre as áreas de aplicações temos vários campos da ciência e da engenharia, incluindo escoamento de fluidos, meios viscosos, reologia, transporte difusivo, redes elétricas, teoria eletromagnética, teoria do campo, probabilidade, etc \cite{mainardi, Hilfer, SATM}. Neste sentido, a construção de métodos eficientes para a solução de equações diferenciais, principalmente os casos não-lineares é um grande desafio para modelagem matemática, em particular quando trabalhamos com o Cálculo Fracionário, pois ainda existem poucas teorias sobre as soluções numéricas e analíticas.  Além disso, do ponto de vista prático, o carácter não local desses operadores impõe dificuldades técnicas para a obtenção de resultados numéricos de modelos descrevendo sistemas realistas. Por exemplo, para obter soluções numéricas de equações diferenciais parciais, frequentemente limita-se o modelo matemático impondo-se condições de fronteira simples (por exemplo, retangular) que não descrevem fielmente o sistema real \cite{N1,N2,N3}. Outro exemplo de dificuldade técnica surge do fato de que existem diversas definições para derivadas fracionárias, e para cada definição, é necessário aplicar um método numérico adequado \cite{WEILBEER,diethelm2005algorithms}. 
  
  Com o objetivo de contornar essas dificuldades, e obter um método numérico que possa ser empregado na resolução de equações diferenciais fracionárias envolvendo diversas formulações de derivadas, no presente trabalho ocorre uma construção de um algoritmo baseado no Método de Monte Carlo para a resolução de equações diferenciais ordinárias fracionárias. O uso do Método de Monte Carlo na resolução de equações diferenciais ordinárias usuais (com derivadas de ordem inteira) é um tema pouco estudado na literatura \cite{MC2} e ainda não explorado no contexto do Cálculo Fracionário. O algoritmo proposto neste trabalho para equações diferenciais fracionárias envolvendo derivadas de Caputo pode ser aplicado em problemas envolvendo outras formulações de derivadas fracionárias e para sistemas de equações diferenciais. 
  

\section{Cálculo Fracionário}
Uma breve introdução à teoria matemática do Cálculo Fracionário será apresentada tendo a definição de integral fracionária proposta por Riemman-Liouville (RL) e as definições mais famosas para o Cálculo Fracionário que serão utilizadas neste artigo.

\begin{defTEMAp}[Integral fracionária de Riemman-Liouville]
Seja $\alpha \in \mathbb{R}_+$, o operador $_aJ^{\alpha}_x$ definido em $L_1([a,b])$ por
\begin{equation}\label{integralRL}
\label{a3}
_a J_x^{\alpha}f(x)= \frac{1}{\Gamma(\alpha)}\int_{a}^{x}(x-u)^{\alpha-1}f(u)du,
\end{equation}
é chamado de integral fracionária de Riemann-Liouville à esquerda, onde $\Gamma(\cdot)$ é a função Gamma, e $a,b\in \mathbb{R}$ com $a<b$.
\end{defTEMAp}
Para $\alpha=n$ inteiro, a integral fracionária de RL (\ref{a3}) coincide com a integral usual de Riemman repetida $n$ vezes \cite{KD}. Da definição (\ref{a3}) pode-se observar que a integral de RL existe para qualquer função $f$ integrável se $\alpha>1$. Além disso, é possível provar a existência da integral de RL  (\ref{a3}) para $f\in L_1([a,b])$ mesmo quando $0<\alpha<1$ \cite{KD}.

A integral fracionária de RL (\ref{a3}) tem um papel central na definição das derivadas fracionárias de RL e Caputo. Para definir as derivadas de RL, onde para inteiros positivos $n>m$, vale a identidade $D^m_x f(x)=D^{n}_x {_aJ}^{n-m}_x f(x)$, onde $D^m_x$ é uma derivada $\frac{d^m}{dx^m}$ de ordem $m$.
\begin{defTEMAp}[Derivada de Riemann-Liouville]
A derivada fracionária à esquerda de Riemann-Liouville de ordem $\alpha \in \mathbb{R}_+$ é dfinida por ${_aD}^{\alpha}_x f(x) = D^{n}_x {_aJ}^{n-\alpha}_x f(x)$, onde $n=[\alpha]+1$, ou seja:
\begin{equation}
\label{a5}
{_aD}^{\alpha}_x f(x)
=\frac{1}{\Gamma(n-\alpha)}\frac{d^n}{dx^n}
\int_{a}^x \frac{f(u)}{(x-u)^{1+\alpha-n}}du,
\end{equation}
onde $a \le x \le b$.
\end{defTEMAp}

Por outro lado, a derivada fracionária de Caputo é definida invertendo-se a ordem entre derivadas e integrais.
\begin{defTEMAp}[Derivada de Caputo]
A derivada fracionária de Caputo à esquerda de ordem $\alpha\in \mathbb{R}_+$ é definida por ${_a^C D}^{\alpha}_x f(x) := {_aJ}^{n-\alpha}_x D^{n}_x f(x)$ com $n=[\alpha]+1$, ou seja:
\begin{equation}
\label{a7}
{_a^C D}^{\alpha}_x f(x) := \frac{1}{\Gamma(n-\alpha)}
\int_{a}^x \frac{f^{(n)}(u)}{(x-u)^{1+\alpha-n}}du,
\end{equation}
onde $a \le x \le b$ e $f^{(n)}(u)=\frac{d^n f(u)}{du^n} \in L_1([a,b])$.
\end{defTEMAp}

Uma consequência importante das definições (\ref{a5}) e (\ref{a7}) é que as derivadas fracionárias de RL e Caputo são operadores não-locais. Essas derivadas à esquerda dependem dos valores da função à esquerda de $x$, isto é, $a\leq u \leq x$. Por outro lado, é importante notar que quando $\alpha=n$ é um inteiro, as derivadas fracionárias de RL e Caputo de reduzem à derivadas de ordem inteira $n$ \cite{KD}. Uma outra consequência dessas definições é que as derivadas fracionárias de RL e Caputo são conectadas pela seguinte relação:
\begin{equation}
\label{a8a}
{_a^C D}^{\alpha}_x f(x)={_a D}^{\alpha}_x f(x)-\sum_{k=0}^{n-1}\frac{f^{(k)}(a)}{\Gamma(k-\alpha+1)}(x-a)^{k-\alpha},
\end{equation}
onde $n=[a]+1$.

Finalmente, a integral e a derivada fracionárias de RL satisfazem a seguinte generalização fracionária do Teorema Fundamental do Cálculo:
\begin{teoTEMA}[Teorema Fundamental do Cálculo de Riemann-Liouville]
\label{thm:1}
Seja $0<\alpha<1$, e seja ${_a J}^{1-\alpha}_x f(x)$ absolutamente contínua em $[a,b]$. A seguinte igualdade é satisfeita:
\begin{equation}
\label{a9}
{_a J}^{\alpha}_b {_a D}^{\alpha}_x f(x) = f(b)-\frac{(x-a)^{\alpha-1}}{\Gamma(\alpha)}\lim_{u\rightarrow a^+}{_a J}^{1-\alpha}_u f(u).
\end{equation}
\end{teoTEMA}
Por outro lado, a generalização do Teorema Fundamental do Cálculo para Caputo é dada por:
\begin{teoTEMA}[Teorema Fundamental do Cálculo de Caputo]
\label{thm:2}
Seja $0<\alpha<1$, e $f$ uma função diferenciável em $[a,b]$. A seguinte igualdade é satisfeita:
\begin{equation}
\label{a11}
{_a J}^{\alpha}_b {_a^C D}^{\alpha}_x f(x) = f(b)-f(a).
\end{equation}
\end{teoTEMA}
A prova dos Teoremas \ref{thm:1} e \ref{thm:2} pode ser vista, por exemplo, em \cite{KD}.
\section{O Método de Monte Carlo para a solução de equações diferenciais ordinárias fracionárias}
O Método de Monte Carlo (MMC) é uma ferramenta poderosa na resolução de diversos problemas, especialmente em física, biologia, engenharias e finanças \cite{MC}. Sua origem remonta ao trabalho pioneiro de Comte de Buffon que, em 1777, propôs um método estatístico para calcular o valor do número $\pi$ baseado no lançamento de uma agulha \cite{MC}. O MMC foi popularizado durante a Segunda Guerra Mundial quando foi utilizado para o projeto da bomba atômica \cite{MC}. Atualmente, o MMC representa uma grande classe de métodos estatísticos que utilizam amostras de números aleatórios para obter a solução numérica de diversos tipos de problemas.

\begin{figure}[!h]
\centering
\subfigure[1 ponto. \label{sim11}]{
\includegraphics[width=.4\textwidth]{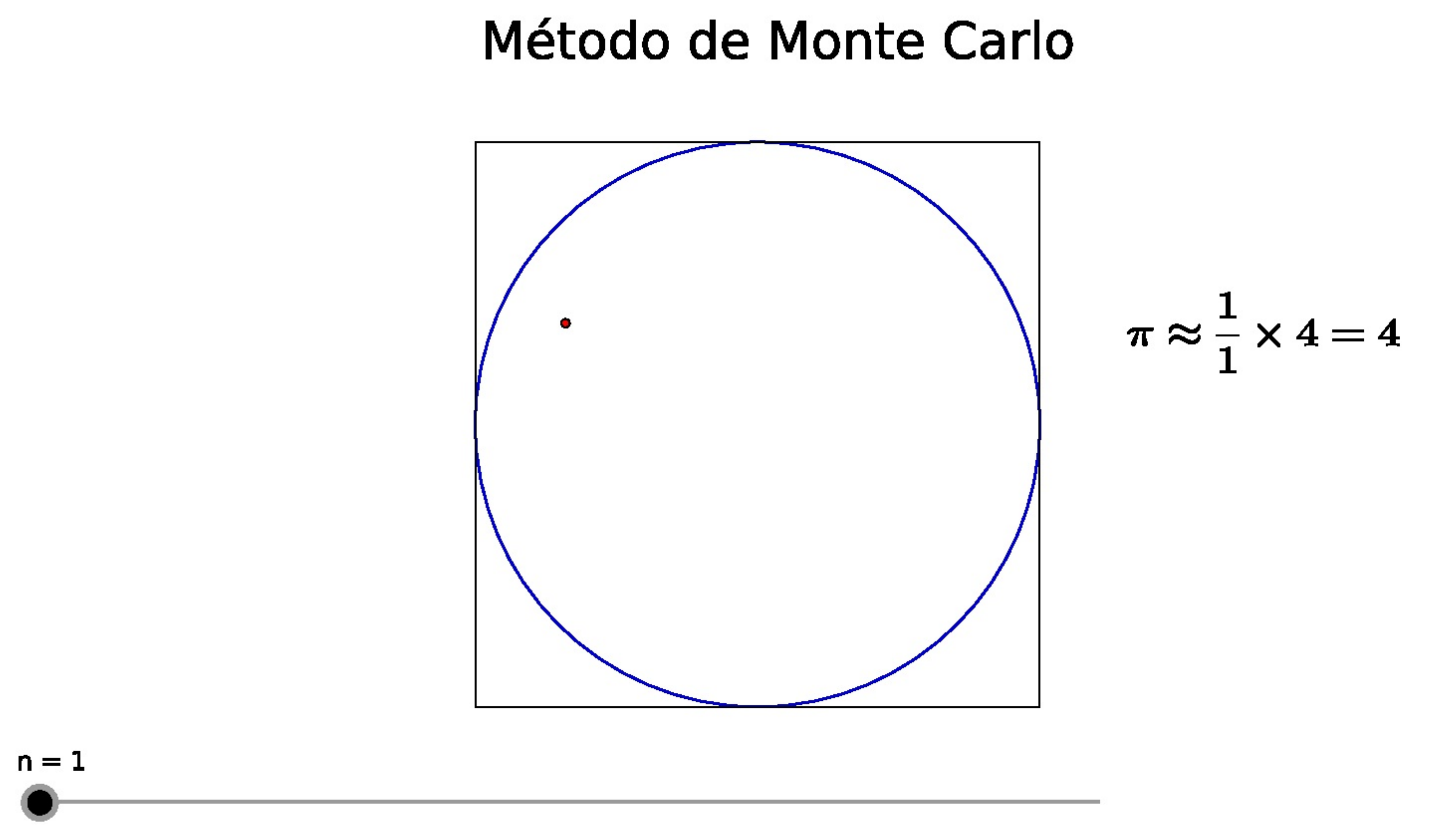}
} 
\qquad 
\subfigure[50 pontos. \label{sim12}]{
\includegraphics[width=.4\textwidth]{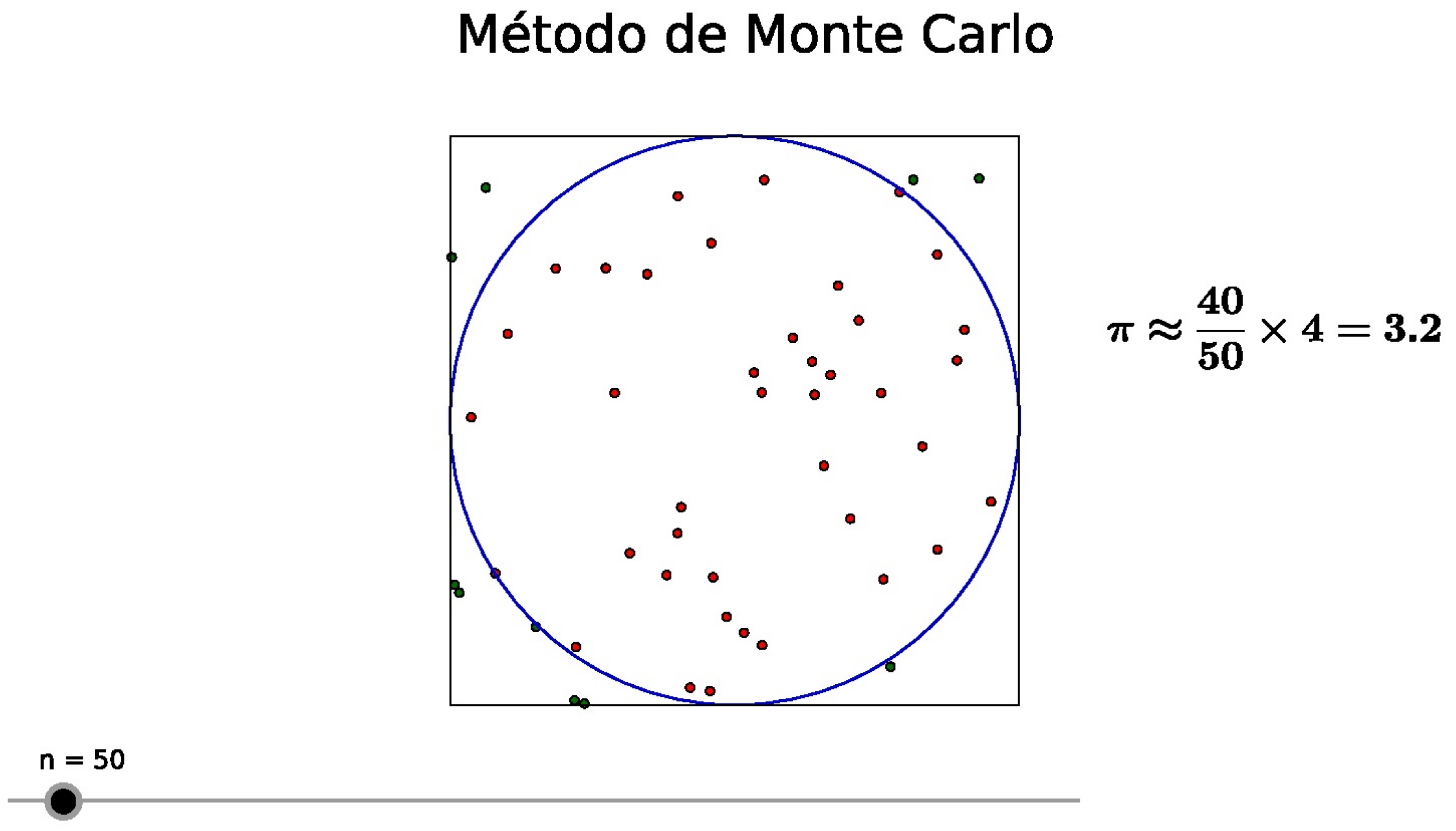}
} 
\caption{Simulação cálculo do $\pi$ com 1 e 50 pontos.}
\label{sim1}
\end{figure}

Uma das aplicações mais simples e interessantes do MMC é o cálculo do $\pi$. Para isso, desenhamos um círculo de raio $1$ dentro de um quadrado de lado $2$ (veja a Figura \ref{sim1}). Ao sortear um ponto aleatório no interior do quadrado, e uniformemente distribuído, a probabilidade $P$ desse ponto estar no interior do círculo é $P=\frac{A_{cir}}{A_{quad}}=\frac{\pi}{4}$, onde $A_{cir}=\pi$ é a área do círculo e $A_{quad}=4$ é a área do quadrado. Portanto, sorteando aleatoriamente $N$ pontos, e sabendo que desses pontos $N_{acertos}$ estão no interior do círculo, podemos obter o valor aproximado de $\pi$ como:
\begin{equation}
\pi \approx \frac{N_{acertos}}{N}A_{quad}=\frac{N_{acertos}}{N}4.
\label{exemplo}
\end{equation}

\begin{figure}[!h]
\centering
\subfigure[100 pontos. \label{sim21}]{
\includegraphics[width=.4\textwidth]{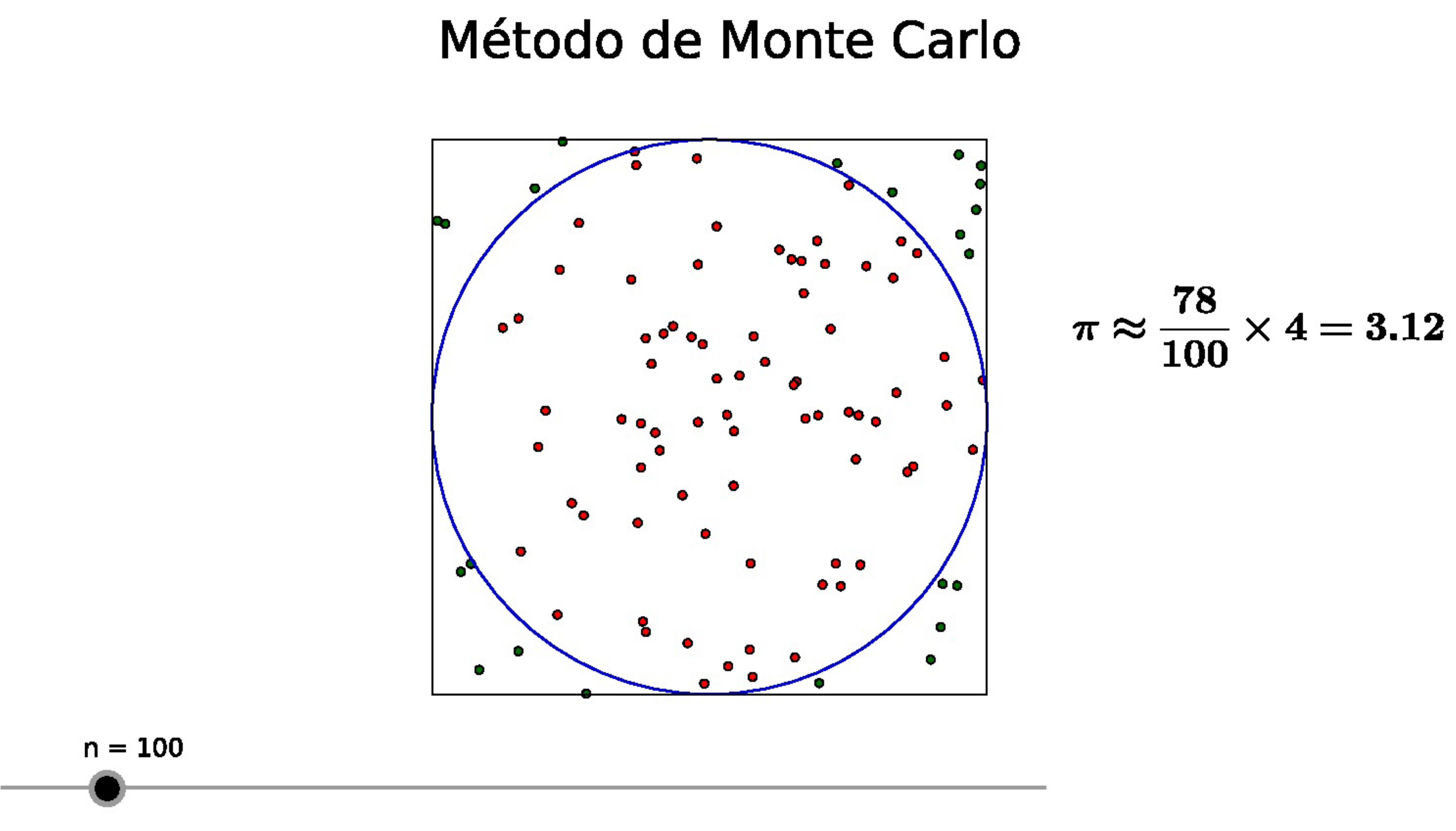}
} 
\qquad 
\subfigure[1000 pontos. \label{sim22}]{
\includegraphics[width=.4\textwidth]{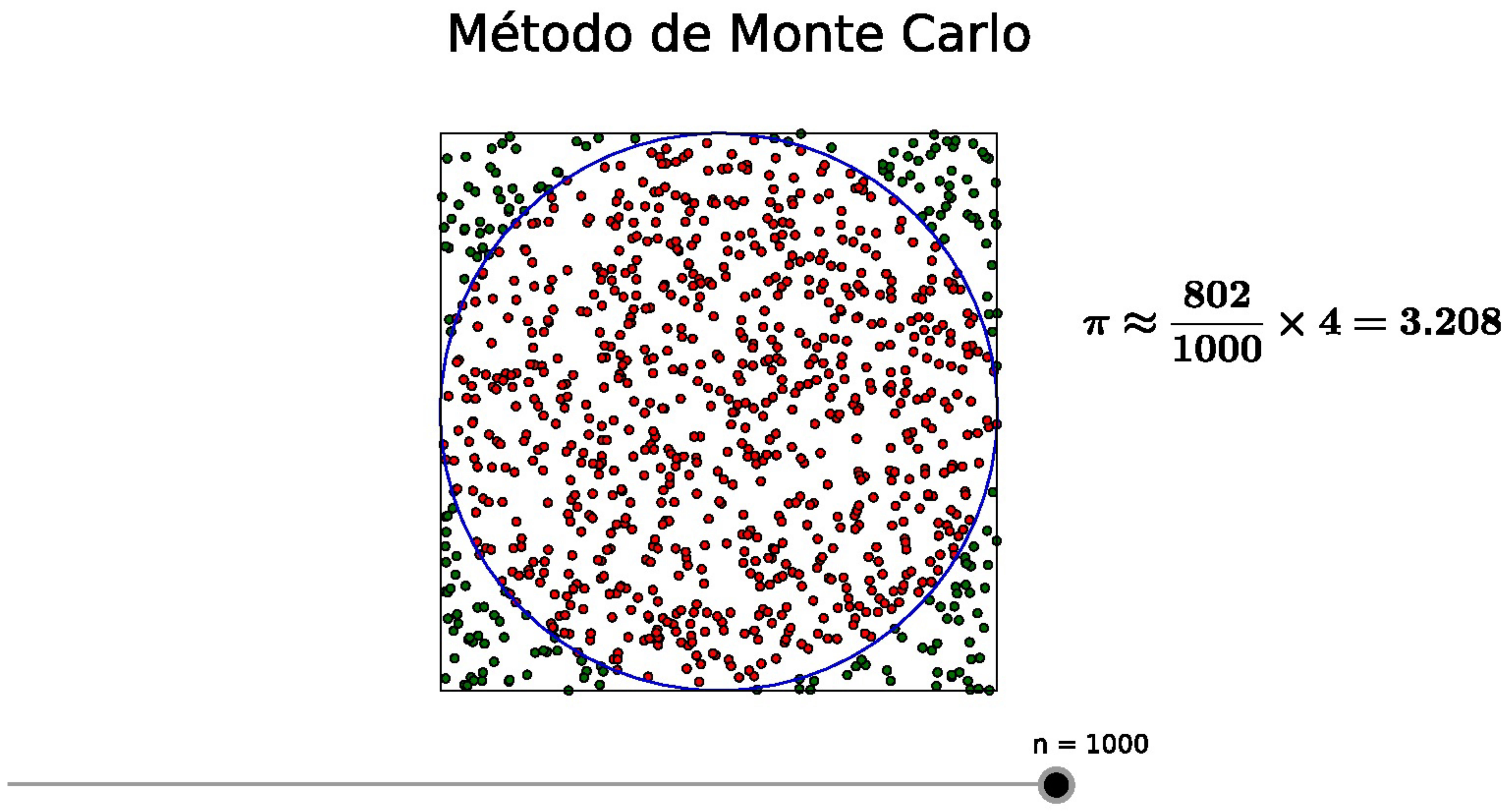}
} 
\caption{Simulação cálculo do $\pi$ com 100 e 1000 pontos.}
\label{sim2}
\end{figure}

 Nas figuras \ref{sim1} e \ref{sim2}, o Método de Monte Carlo é implementado junto com a equação (\ref{exemplo}) para obter valores aproximados para a constante $\pi$, utilizando 1, 50, 100 e 1000 pontos sorteados aleatoriamente, e com distribuição uniforme. 

Uma outra aplicação importante é o cálculo de integrais definidas. Um MMC simples para calcular integrais envolve o uso do Teorema do Valor Médio para integrais. Seja $f(x)$ uma função integrável no sentido de Riemman no intervalo $[a,b]$. Então
\begin{equation}
\label{b1}
\int_{a}^{b}f(x)dx=(b-a)\langle f\rangle,
\end{equation}
onde $\langle f\rangle$ é o valor médio da função em $[a,b]$. Neste ponto, uma amostra de números aleatórios pode ser utilizada para calcular aproximadamente $\langle f\rangle$. Dado um conjunto de $N$ números ($x_1,x_2,...,x_N$) aleatórios uniformemente distribuídos em $[a,b]$, um valor aproximado $\langle f\rangle_N$ para $\langle f\rangle$ é obtido simplesmente calculando $f(x)$ para cada número aleatório selecionado. Nós temos:
\begin{equation}
\label{b2}
\langle f\rangle_N = \frac{1}{N}\sum_{i=1}^N f(x_i),
\end{equation}
onde $\langle f\rangle_N$ converge estatisticamente para $\langle f\rangle$ no limite $N\rightarrow \infty$ \cite{MCb,MC3}, ou seja, $\langle f\rangle = \lim_{N\rightarrow \infty}\langle f\rangle_N$. Portanto, (\ref{b2}) nos fornece uma aproximação para a integral (\ref{b1}):
\begin{equation}
\label{b3}
\int_{a}^{b}f(x)dx \approx \frac{b-a}{N}\sum_{i=1}^N f(x_i),
\end{equation}
que estatisticamente converge para o valor exato da integral quando $N\rightarrow \infty$. Como a variância $\sigma^2$ para $f$ é dada por \cite{MCb,MC3}:
\begin{equation}
\label{b4}
\sigma^2=\frac{b-a}{N}\sum_{i=1}^N f^2(x_i)-\left(\frac{b-a}{N}\sum_{i=1}^N f(x_i)\right)^2,
\end{equation}
e o erro padrão $Err$ é:
\begin{equation}
\label{b5}
Err=\frac{\sqrt{\sigma^2}}{\sqrt{N}},
\end{equation}
a aproximação (\ref{b3}) dá um resultado para a integral que difere por até um erro padrão do valor exato com probabilidade de $68.3\%$ (e difere por até dois erros padrões do valor exato com probabilidade de $95.4\%$). Como o erro estimado é proporcional à variância e inversamente proporcional à raiz quadrada de $N$, existem dois métodos de reduzir o erro. O primeiro é aumento $N$, e o segundo é reduzindo a variância usando $N$ números aleatórios ($x_1,x_2,...,x_N$) não-uniformemente distribuídos em $[a,b]$ \cite{MCb,MC3}.

\subsection{Resolução de equações diferenciais fracionárias}

Neste trabalho um MMC é proposto para obter a solução numérica de equações diferenciais ordinárias fracionárias com derivadas de Caputo. Considerando o seguinte problema de valor inicial: 
\begin{equation}
\label{edo1}
_a^{C} D^{\alpha}_x y = f(x, y),\;\;\; y(a)=y_a,
\end{equation}
com $0<\alpha \leq 1$. O primeiro passo para obter a solução do problema de valor inicial fracionário (\ref{edo1}) consiste no uso do Teorema Fundamental do Cálculo de Caputo. Da equação (\ref{a11}) temos que ${_a J}^{\alpha}_x {_a^C D}^{\alpha}_x y(x)=y(x)-y(a)=y(x)-y_a$. Portanto, assumindo que $y(x)$ é diferenciável no domínio, e integrando os dois lados da (\ref{edo1}) com uma integral fracionária de RL (\ref{a3}), temos:
\begin{equation}
\label{edo2}
y(x)=y_a+\frac{1}{\Gamma(\alpha)}\int_{a}^{x}(x-u)^{\alpha-1}f(u,y(u))du.
\end{equation}
O segundo passo consiste no uso do MMC para obter a solução da integral em (\ref{edo2}). No entanto, como $y(x)$ não é conhecida, não se pode aplicar o MMC diretamente para o cálculo dessa integral. Define-se primeiro uma discretização da função $y(x)$. Para um inteiro positivo $L$, temos $x_n=x_0+nh$ ($n=0,1,...,L$), onde $h=\frac{x_L-x_0}{L}$ e $x_0=a$. Seja então $y_n=y(x_n)$, a discretização da função $y(x)$ é definida por:
\begin{equation}
\label{edo3}
y^L(x)=y_n\;\;\; \mbox{se} \;\;\; x_{n}\leq x < x_{n+1}.
\end{equation}
Como foi assumido que $y(x)$ é diferenciável, temos $\lim_{L\rightarrow \infty} y^L(x)=y(x)$. Portanto, a função discretizada $y^L(x)$ é uma boa aproximação da função $y(x)$ para $L$ suficientemente grande. Pode-se então encontrar os valores $y_n$ à partir da seguinte relação de recorrência:
\begin{equation}
\label{edo4}
y_{n}=y_a+\frac{1}{\Gamma(\alpha)}\int_{a}^{x_n}(x_n-u)^{\alpha-1}f(u,y^L(u))du,\;\;\; (n=1,...,L).
\end{equation}
Finalmente, calculando a integral em (\ref{edo4}) usando o MMC (como feito na (\ref{b3})), obtêm-se a solução aproximada do problema de valor inicial (\ref{edo1}):
\begin{equation}
\label{edo6}
y_{n}\approx y_a+\frac{1}{\Gamma(\alpha)}\frac{x_n-a}{N}\sum_{i=1}^{N}(x_n-u_i)^{\alpha-1}f(u_i,y^L(u_i)),
\end{equation}
onde $N$ é um inteiro positivo e $(u_1,...,u_N)$ são números aleatórios uniformemente distribuídos em $[a,x_n]$. Finalmente, é importante ressaltar que o MMC proposto pode ser utilizado também para obter a solução de equações diferenciais ordinárias com derivada inteira fazendo-se $\alpha=1$ na (\ref{edo6}).
\subsection{Resultados numéricos}

Será apresentado agora alguns exemplos de solução de problemas de valores iniciais fracionários utilizando o MMC proposto. Nestes exemplos, utilizou-se $L=100$ e $N=10$ na equação (\ref{edo6}) para o cálculo do valor aproximado de $y_n$. Por questões didáticas, para destacar o caráter estatístico do método, optou-se por um valor de $N$ pequeno para que as divergências estatísticas entre a solução numérica e a solução analítica fiquem graficamente visíveis nas figuras. O método proposto pode ser aplicado para problemas de valores iniciais (\ref{edo1}) gerais, mesmo guando $f(x,y)$ é não linear em $y$. No entanto, os exemplos apresentados são dos casos mais simples onde podemos comparar a solução numérica com a solução analítica do problema.

\subsubsection{Exemplo 1}
\begin{figure}[!h]
    \centering
    \includegraphics[width=0.9\textwidth]{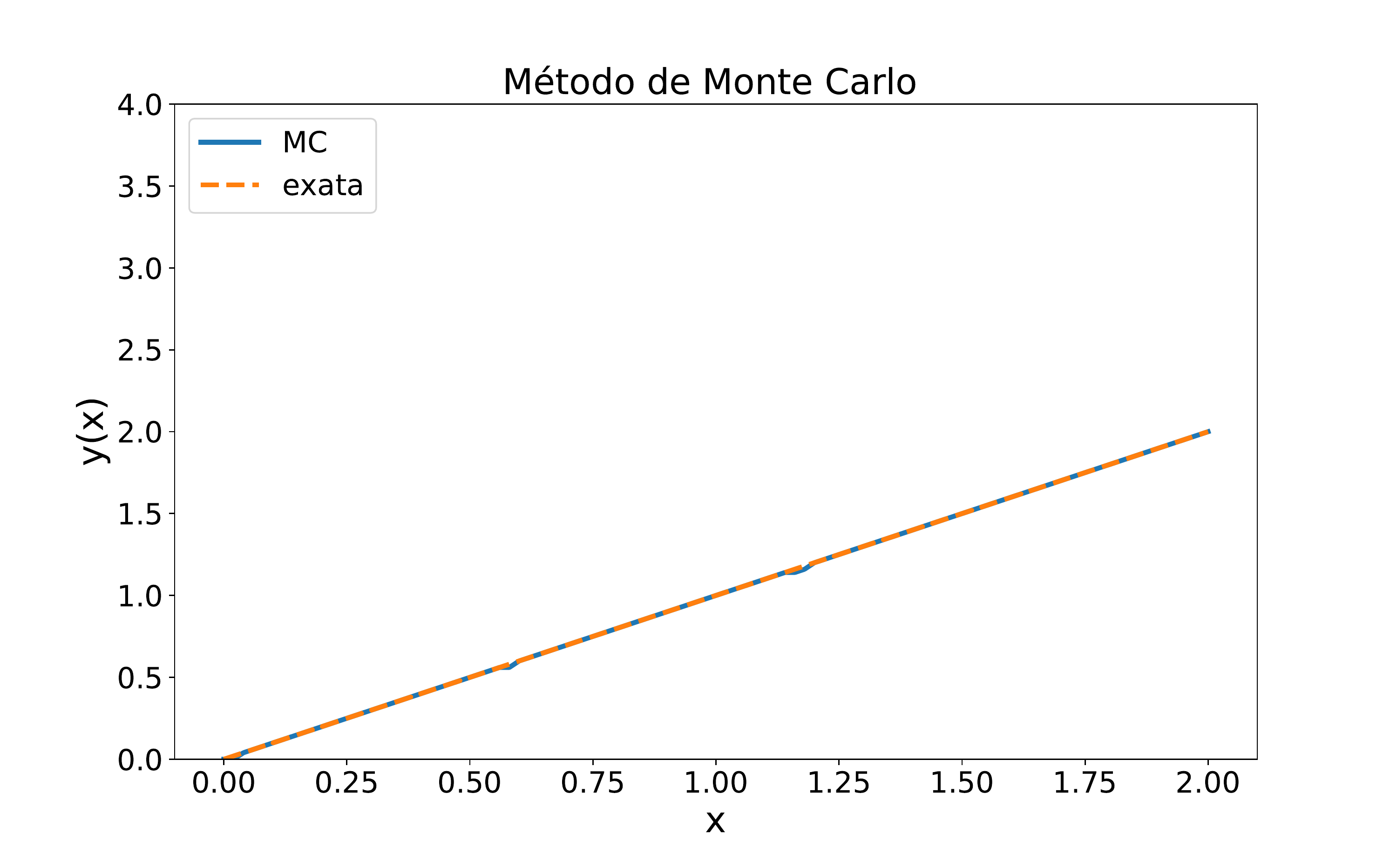}
    \caption{Comparação da solução exata com a obtida pelo MMC para $\alpha=1$ e $\beta=0$ da equação \ref{exemplo1}. \\ Fonte: do autor.} 
    \label{f01}
\end{figure}

\begin{figure}[!htpb]
    \centering
    \includegraphics[width=0.9\textwidth]{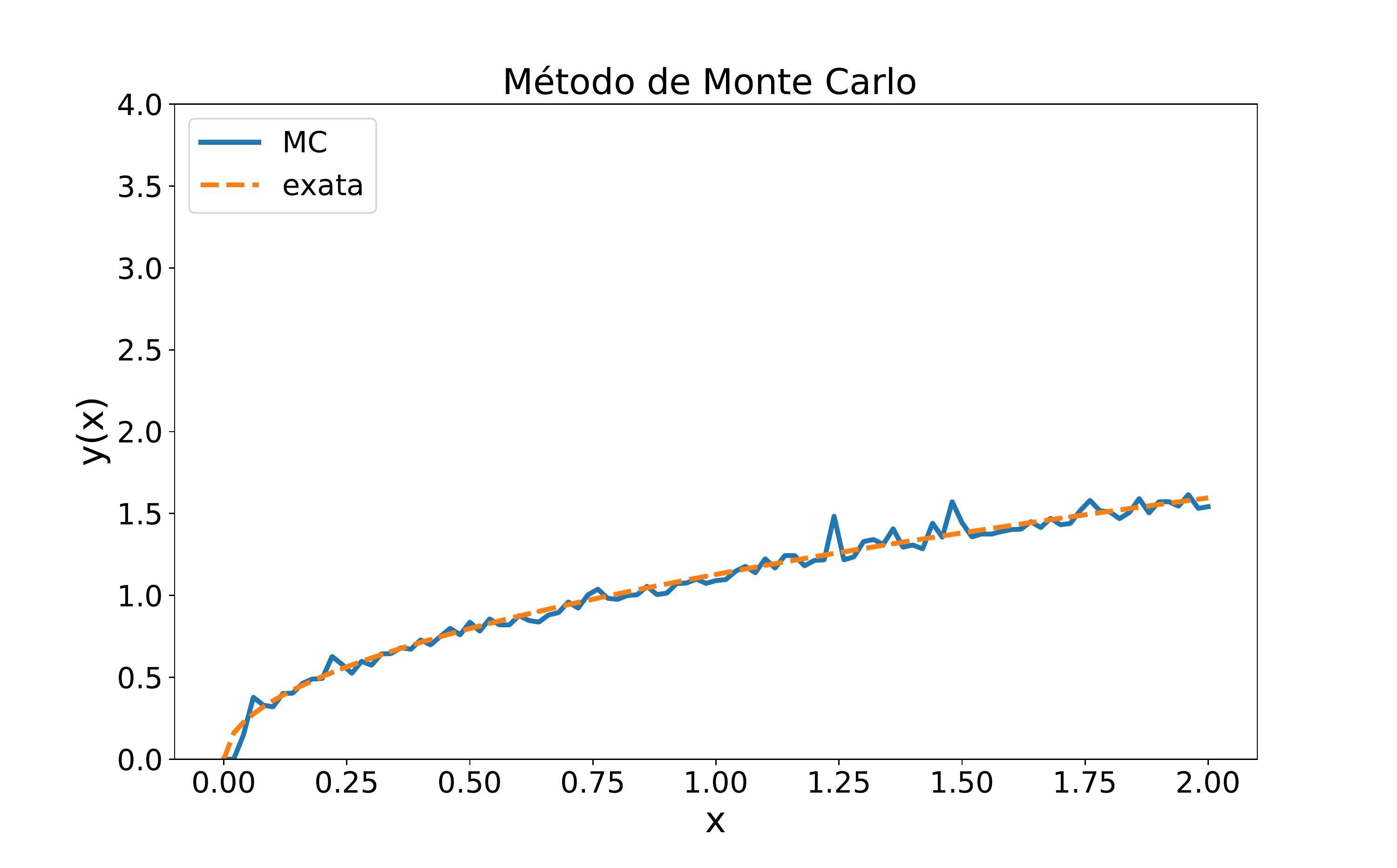}
    \caption{Comparação da solução exata com a obtida pelo MMC para $\alpha=0.5$  $\beta=0$ da equação \ref{exemplo1}.\\ Fonte: do autor.} 
    \label{f02}
\end{figure}

\begin{figure}[!htpb]
    \centering
    \includegraphics[width=0.9\textwidth]{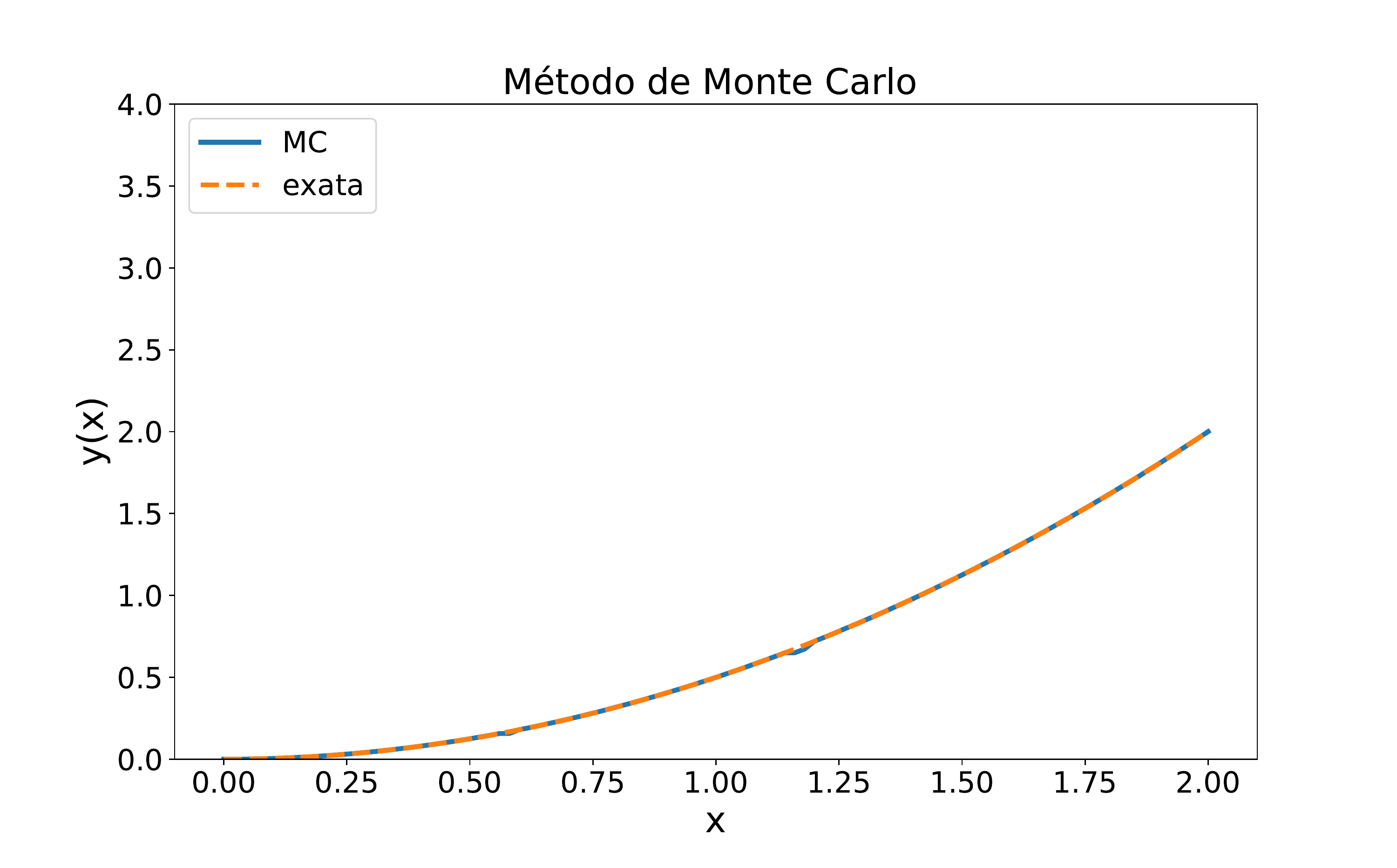}
    \caption{Comparação da solução exata com a obtida pelo MMC para $\alpha=1$ e $\beta=1$ da equação \ref{exemplo1}. \\ Fonte: do autor.} 
    \label{f03}
\end{figure}

\begin{figure}[!htpb]
    \centering
    \includegraphics[width=0.9\textwidth]{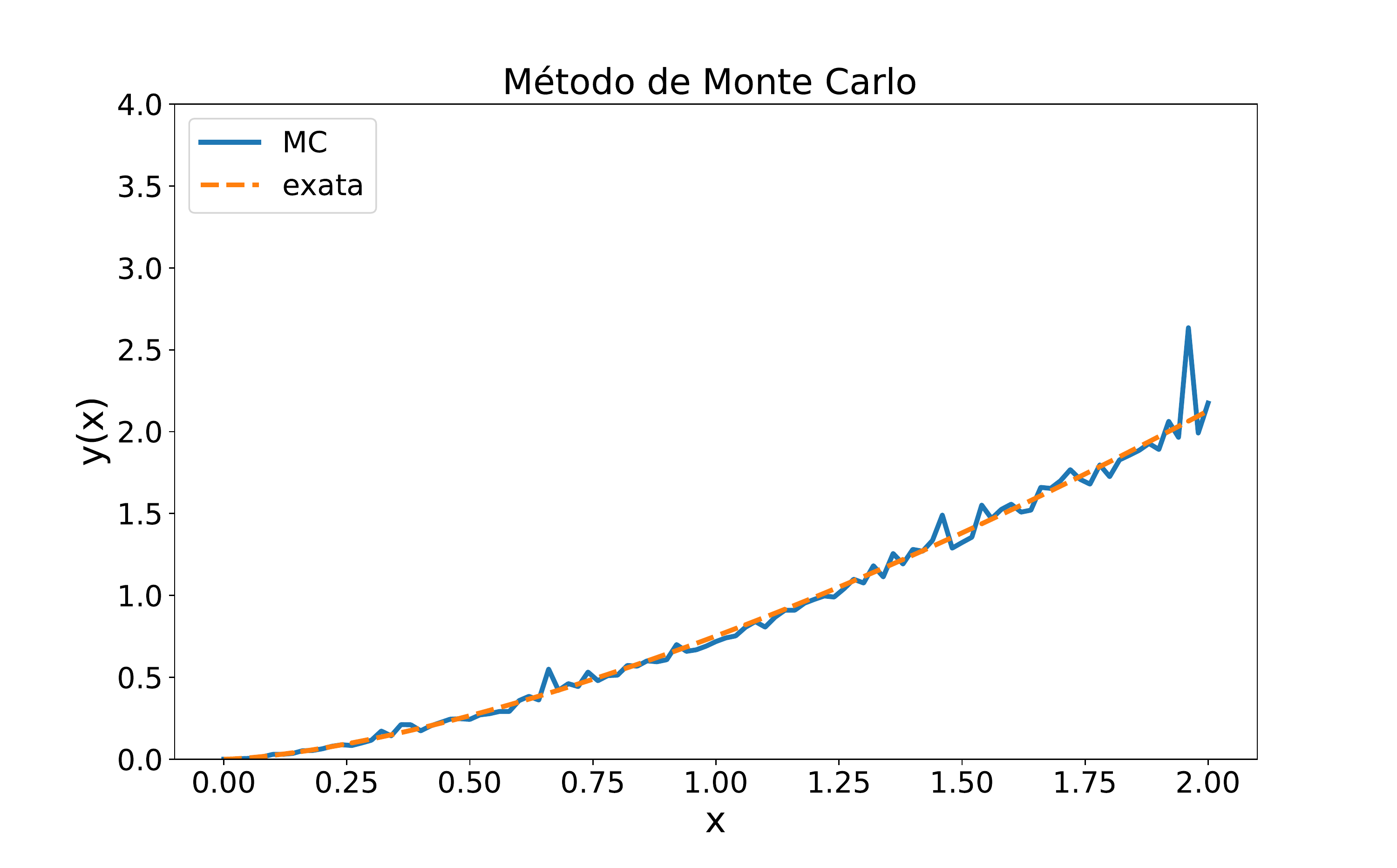}
    \caption{Comparação da solução exata com a obtida pelo MMC para $\alpha=0.5$ e $\beta=1$ da equação \ref{exemplo1}. \\ Fonte: do autor.} 
    \label{f04}
\end{figure}

O primeiro exemplo é o Problema de Valor Inicial (PVI)
\begin{equation}
\label{exemplo1}
_a^{C} D^{\alpha}_x y=x^{\beta},\;\;\; y(0)=0
\end{equation}
com $0<\alpha\leq 1$. A solução analítica do problema, que pode ser obtida usando transformada de Laplace, é dada por:
\begin{equation}
\label{solexemplo2}
y(x)= \frac{x^{\alpha +\beta}}{\Gamma(\alpha + \beta+1)}.
\end{equation}
Nas figuras \ref{f01},\ref{f02},\ref{f03} e \ref{f04}, compara-se a solução exata (\ref{solexemplo2}) com a solução obtida pelo MMC proposto para diferentes valores de $\alpha$ e de $\beta$. 

Na Figura \ref{f01} vemos uma boa concordância entre a solução exata e a solução pelo MMC para o caso usual de derivada inteira ($\alpha=1$) e $\beta=1$. Como discutido na seção anterior, o MMC proposto pode também ser aplicado para a solução de equações diferenciais ordinárias de ordem inteira. As pequenas discrepâncias entre a solução exata e a solução pelo MMC podem ser reduzidas aumentando-se o valor de $N$.
A solução para o caso de derivada fracionária de ordem $\alpha=0.5$, e $\beta=0$, é apresentada na Figura \ref{f02}. Neste caso de derivada fracionária, existe também uma boa concordância entre a solução exata e a solução obtida pelo MMC, apesar de uma maior discrepâncias (devido a uma pior convergência estatística) entre a solução exata e a solução pelo MMC. Essa discrepância maior é devida ao polo da função $(x_n-u_i)^{\alpha-1}$ em (\ref{edo6}). As alternativas para se reduzir o erro, é aumentar o valor de $N$ ou utilizar distribuições de probabilidades não-uniformes \cite{MCb,MC3}.

Nas Figuras \ref{f03} e \ref{f04}, para $\beta=1$, apresenta-se os resultados para derivadas de ordem inteira $\alpha=1$ e ordem fracionária $\alpha=0.5$, respectivamente.
Em ambos os casos, novamente temos uma boa concordância entre s solução exata e a solução obtida pelo MMC. Como no caso anterior, a solução da equação fracionária apresenta maior discrepâncias devido a pior convergência estatística.
\subsubsection{Exemplo 2}
O segundo exemplo é o Problema de Valor Inicial (PVI)
\begin{equation}
\label{a222}
_a^{C} D^{\alpha}_x y=y,\;\;\; y(0)=1
\end{equation}
com $0<\alpha\leq 1$. A solução analítica do problema, que pode ser obtida usando transformada de Laplace, é dada por:
\begin{equation} 
y(x)=\sum_{k=0}^{\infty} \frac{{(x)}^{\alpha k}}{\Gamma(\alpha\,k +1)}=E_{\alpha}(x).
\end{equation}
onde $E_{\alpha}(x)$ é chamada de Função de Mittag-Leffler de um parâmetro \cite{mittag}.
\begin{figure}[!h]
    \centering
    \includegraphics[width=0.9\textwidth]{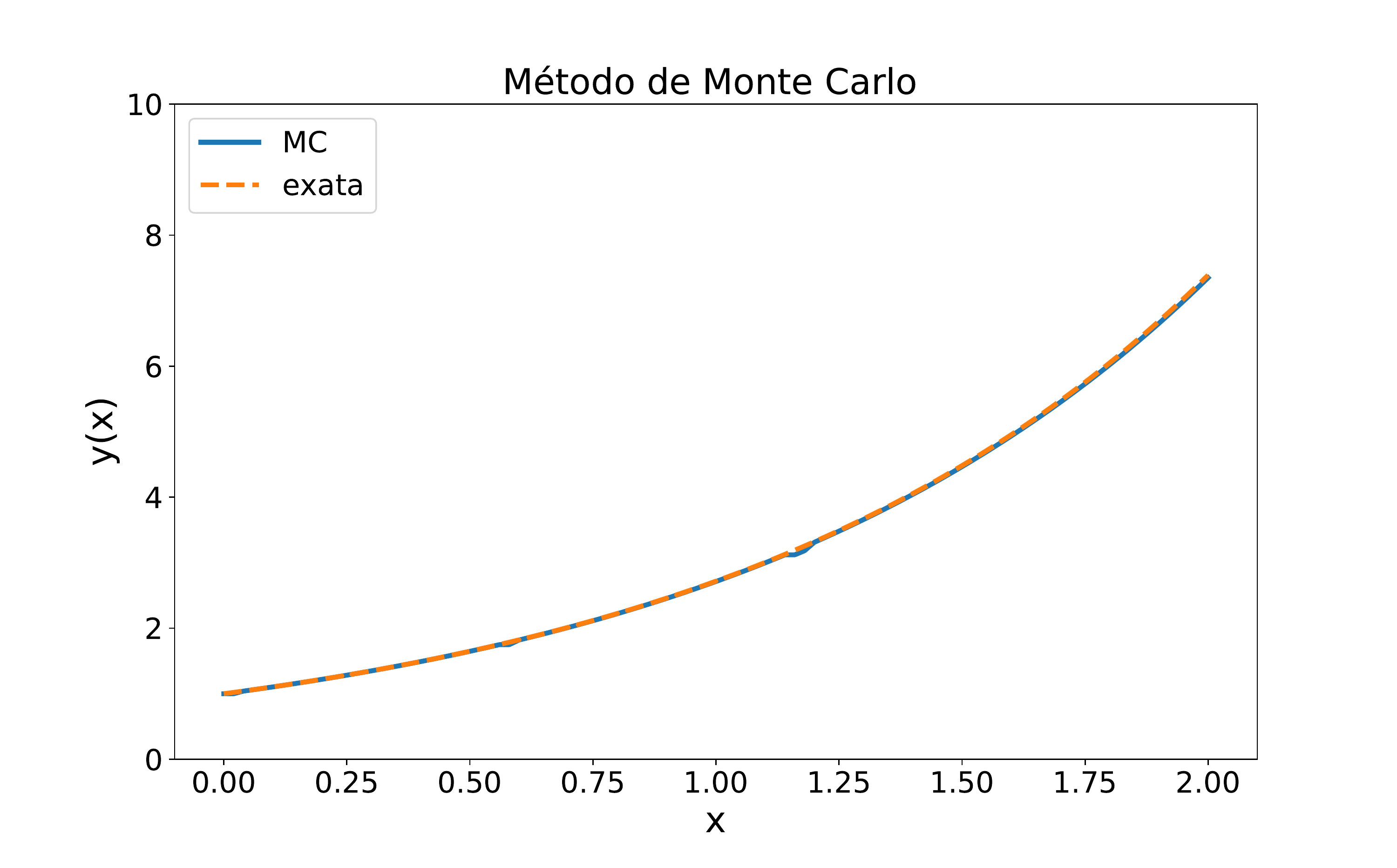}
    \caption{Comparação da solução exata com a obtida pelo MMC para $\alpha=1$ da equação \ref{a222}. \\ Fonte: do autor.} 
    \label{fig:01}
  \end{figure}
  
\begin{figure}[!h]
    \centering
    \includegraphics[width=0.9\textwidth]{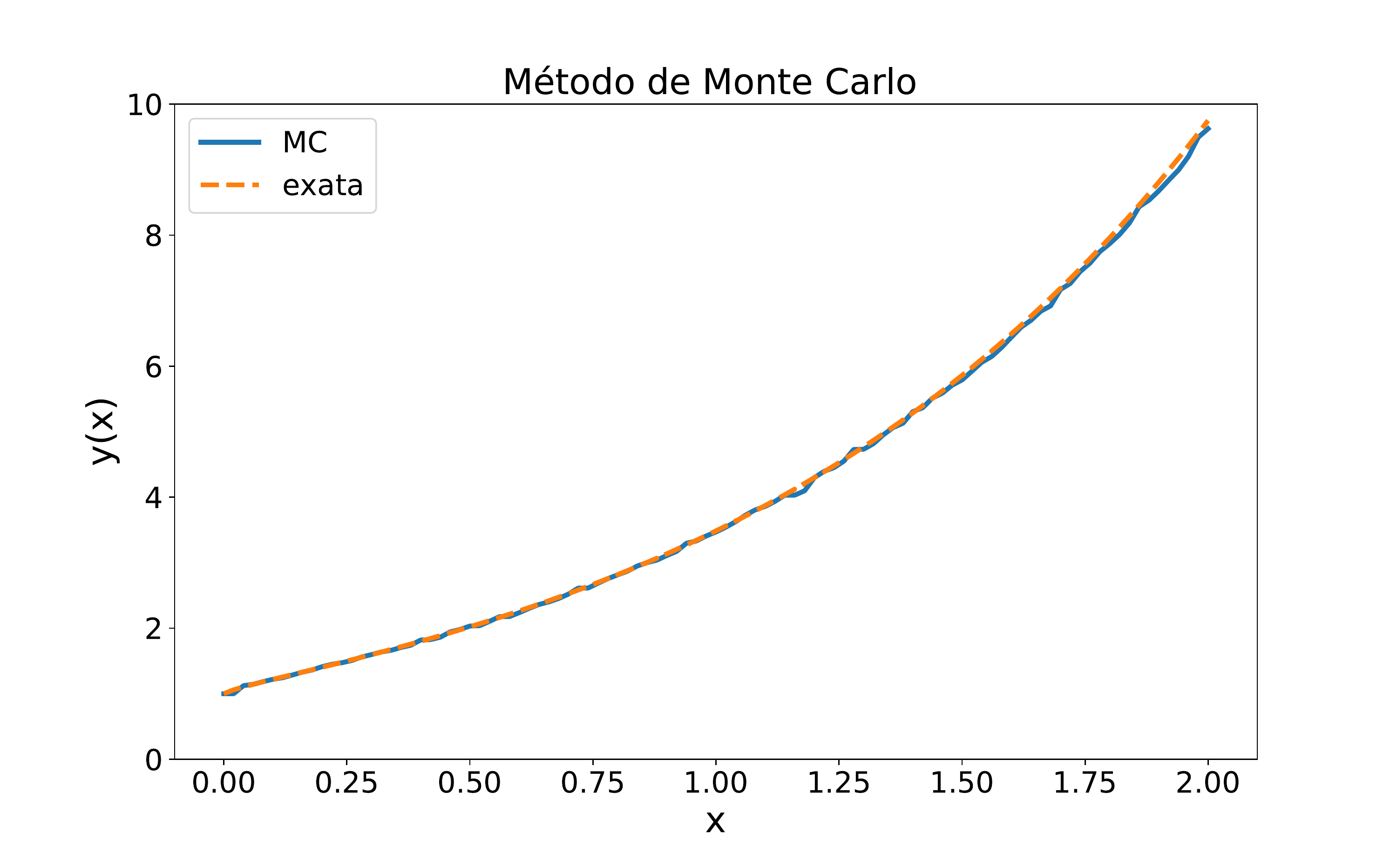}
    \caption{Comparação da solução exata com a obtida pelo MMC para $\alpha=0.75$ da equação \ref{a222}. \\ Fonte: do autor.} 
    \label{fig:02}
  \end{figure}

\begin{figure}[!h]
    \centering
    \includegraphics[width=0.9\textwidth]{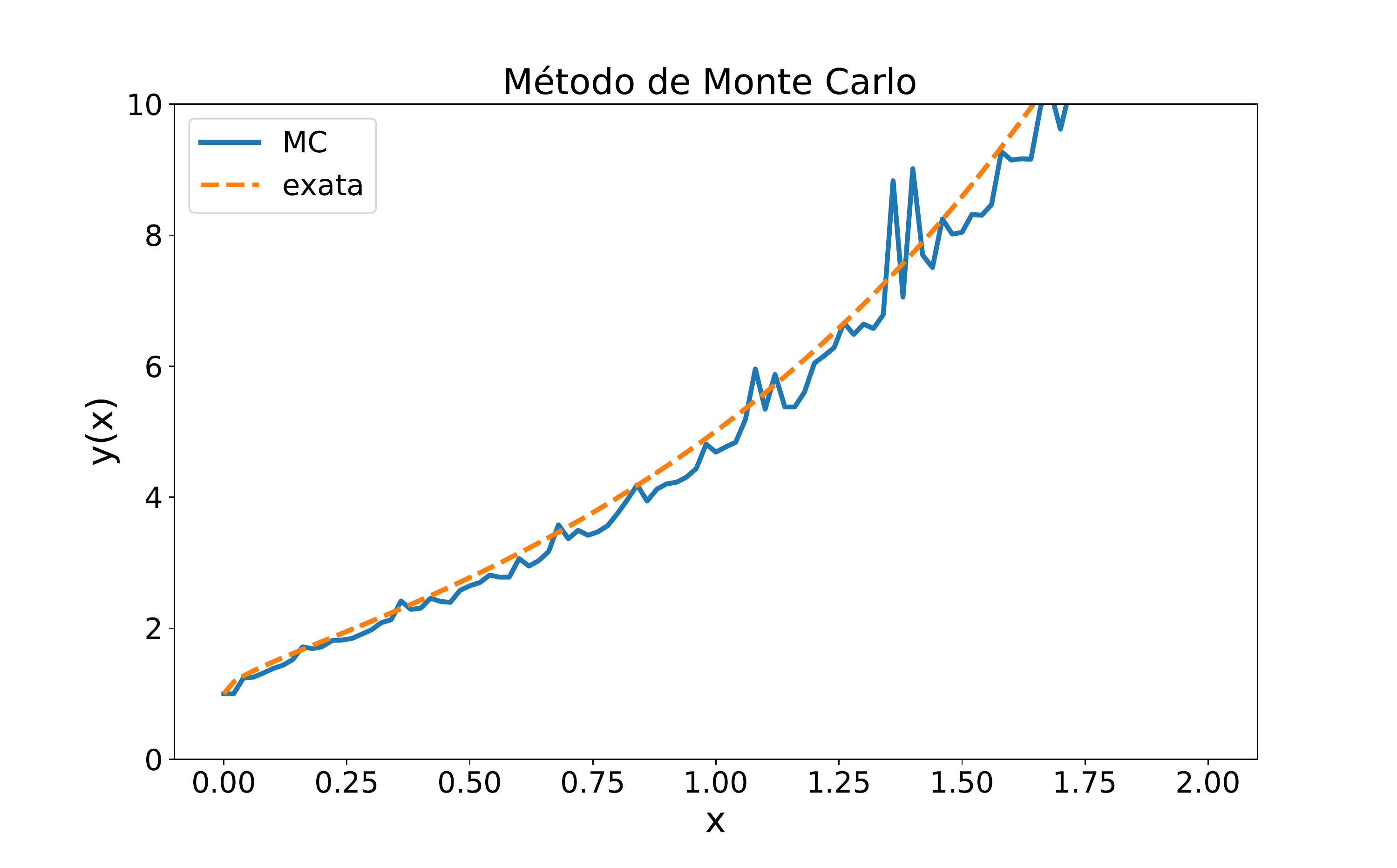}
    \caption{Comparação da solução exata com a obtida pelo MMC para $\alpha=0.5$ da equação \ref{a222}. \\ Fonte: do autor.} 
    \label{fig:03}
  \end{figure}

Na Figura  \ref{fig:01}, temos a solução da equação diferencial com $\alpha=1$ (derivada de ordem inteira), na Figura  \ref{fig:02}, temos a solução de equação diferencial com $\alpha=0.75$ e na Figura \ref{fig:03}, temos a solução da equação diferencial com $\alpha=0.5$. Nos gráficos confrontamos os valores exatos e os obtidos pelo MMC. 
  Em todas essas figuras, observa-se uma boa concordância entre a solução exata e a solução pelo MMC. Observa-se ainda, que quanto menor a ordem da derivada (menor o valor de $\alpha$), pior a convergência estatística do método. Para obter as soluções apresentadas nos gráficos,  o domínio foi divido em $L=100$ partes e foi feito uma amostragem de $10$ pontos aleatórios em cada uma dessas partes ($N=10n$). Esses valores foram escolhidos apenas para ilustrar o funcionamento do método. A convergência estatística da solução de Monte Carlo para a solução exata pode ser melhorada (reduzindo o erro em relação à solução exata) aumentando-se o número de amostragem ou utilizando distribuições não uniformes. Há um estudo em andamento para otimizar esta convergência.

  \section{Conclusões}

  O Cálculo Fracionário é usado com sucesso na modelagem de sistemas complexos \cite{mainardi, Hilfer, SATM, Magin,rivero,random}. Baseado nestes resultados, observar-se que equações diferenciais fracionárias podem ser usadas na modelagem matemática de diversos problemas reais. Neste sentido faz-se necessário a construção de métodos numéricos eficientes para encontrar a solução de equações fracionárias. Um dos métodos que estudamos é o Método de Monte Carlo \cite{MC0, MC1, MC2}. Este método é importante pelo eu carácter universal (pode ser facilmente adaptado para definições diferentes de derivadas fracionárias), e pelo potencial de ser utilizado em problema de fronteira não retangulares. Não foi encontrado na literatura nenhum trabalho que utiliza o Método de Monte Carlo, como o que foi utilizado neste estudo, para resolver equações diferenciais ordinárias fracionárias. Para estudos futuros, pretende-se implementar o método para equações diferenciais parciais fracionárias e estudar o problema de fronteiras não retangulares.

  \section*{Agradecimentos}
  Agradecimentos ao Programa de Pós-Graduação em Modelagem Computacional da Universidade Federal do Rio Grande (PPGMC-FURG) e ao Instituto de Matemática, Física e Estatística da Universidade Federal do Rio Grande (IMEF-FURG) pelo apoio e incentivo. Também agradecemos à FAPERGS e ao CNPq pelo suporte financeiro.
  
\begin{abstract}
{\bf Abstract}. 
This article analyzes and develops a method to solve fractional ordinary differential equations using the Monte Carlo Method. A numerical simulation is performed for some differential equations, comparing the results with what exists in the mathematical literature. The Python language is used to create computational models.
\end{abstract}

\end{document}